\documentclass[reqno,12pt]{amsart}
\usepackage{amsmath,amsthm,amssymb,tabularx}

\theoremstyle{plain}

\newtheorem{tm}{Theorem}[section]
\newtheorem{lm}{Lemma}[section]
\newtheorem{cor}{Corollary}[section]
\theoremstyle{definition}
\newtheorem{ex}{Example}[section]

\newcommand{\beq}{\begin{equation}}
\newcommand{\eeq}{\end{equation}}
\newcommand{\bga}{\begin{gather*}}
\newcommand{\ega}{\end{gather*}}
\newcommand{\bit}{\begin{itemize}}
\newcommand{\eit}{\end{itemize}}
\newcommand{\btm}{\begin{tm}}
\newcommand{\etm}{\end{tm}}
\newcommand{\blm}{\begin{lm}}
\newcommand{\elm}{\end{lm}}
\newcommand{\bcor}{\begin{cor}}
\newcommand{\ecor}{\end{cor}}
\newcommand{\bex}{\begin{ex}}
\newcommand{\eex}{\end{ex}}
\newcommand{\bpr}{\begin{proof}}
\newcommand{\epr}{\end{proof}}

\def\FN{{\mathcal F}(\N)}
\def\N{\mathbb{N}}
\def\R{\mathbb{R}}
\def\Q{\mathbb{Q}}
\def\T{\mathbb{T}}
\def\L{{\mathcal L}}

\def\A{\mathcal{A}}
\def\Z{\mathbb{Z}}
\def\Zp{\Z(p^\infty)}
\let\a\alpha
\let\d\delta
\def\e{\varepsilon}
\let\s\sigma

\def \le {\leqslant}
\def \ge {\geqslant}
\let \<\langle
\let \>\rangle
\let\lf\lfloor
\let\rf\rfloor
\let\phi\varphi

\def\lpw{\L_p^w(G)}

\def\sk#1{\noalign{\kern#1}}

\def\supp{\operatorname{supp}}

\begin{document}

\title
{CONSTRUCTIONS OF REGULAR ALGEBRAS $\L_p^w(G)$}
 \author{Yu. N. Kuznetsova}
 \address{VINITI, Mathematics dept., Usievicha 20, Moscow 125190}
 \email{jkuzn@mccme.ru}
\thanks{Work supported by RFBR grant No. 05-01-00982.}

%\udc{517.986}
\maketitle

\begin{abstract}
Criterion of (Shilov) regularity for weighted algebras $\L_1^w(G)$ on a locally compact
abelian group $G$ is known by works of Beurling (1949) and Domar (1956).
In the present paper this criterion is extended to translation invariant
weighted algebras $\lpw$ with $p>1$. Regular algebras $\lpw$ are constructed
on any sigma-compact abelian group $G$. It was proved earlier by the author
that sigma-compactness is necessary (in the abelian case) for
the existence of weighted algebras $\lpw$ with $p>1$.
\end{abstract}

\section{Criterion of regularity}\label{tauber}

Regular algebras, introduced by G.~E.~Shilov\cite{shilov}, form an important
class of commutative semisimple Banach algebras. Recall their definition.
If a commutative Banach algebra $\A$ is semisimple, the Gel'fand transform
allows to identify it with a subalgebra of continuous function on the
space of maximal ideals $\Sigma$, which is also called the spectrum of $\A$.
Further, $\A$ is said to be regular, if the corresponding algebra of functions
separates points and closed sets in $\Sigma$, i.e. for any closed set
$F\subset \Sigma$ and a point $x\notin F$ there is $f\in\A$ such that
$f(x)\ne0$ and $f$ vanishes on $F$.

Regular algebras posess different interesting properties.  For example,
if we consider the ideal $I(F)$ of all functions in $\A$ that vanish on a
closed subset $F\subset \Sigma$, the set of common zeros of all functions in $I(F)$
equals to $F$ again. If functions with compact support are dense in $\A$ then
every proper closed ideal $J\subset \A$ is contained in a modular maximal
ideal, i.e. has a common zero $s\in\Sigma$ (it is an abstract form of the
Wiener Tauberian theorem, see, e.g., \cite[\S 25D]{loom}).

We start with necessary definitions. $G$ denotes always a locally compact abelian
group, all integrals are taken with respect to Haar measure $\mu$,
$p\ge1$, $1/p+1/q=1$ (for $p=1$ we put $q=\infty$).
A weight is any positive measurable function on $G$.
The space $\L_p^w(G)$ with a weight $w$ is defined as $\{f:fw\in \L_p(G)\}$,
with the norm $\|f\|_{p,w}=\big(\! \int |fw|^p\, \big)^{1/p}$. Indices $p,w$
are sometimes omitted. Weights $w_1$, $w_2$ are called equivalent if
for some constants $C_1$, $C_2$ locally almost everywhere
\beq\label{equiv}
C_1\le {w_1\over w_2} \le C_2.
\eeq
Equivalent weights define the same weighted space and
equi\-va\-lent norms on it.

Recall sufficient conditions on a weight function for $\L_p^w(G)$
to be an algebra with the usual convolution:
$(f*g)(s) = \int f(t)g(s-t)dt$.
For $p=1$ it is submultiplicativity:
\beq\label{eq_submult}
w(s+t) \le w(s)w(t),
\eeq
and for $p>1$ the following inequality (locally almost everywhere):
\beq\label{eqw}
w^{-q} * w^{-q} \le w^{-q}.
\eeq
The space $\lpw$ is translation invariant iff \cite{edw}
for any $s\in G$
\beq\label{ledw}
L_s = {\rm ess}\sup_{t\in G} {w(s+t)\over w(t)}.
\eeq
It is obviously sufficient that $w$ be submultiplicative.

Weighted algebras $\lpw$ are not regular for any weight $w$.
This is related to the fact that for a quickly growing weight $w$ Fourier transforms
of functions in $\L_p^w(\R)$ form a quasianalytic class, for which uniqueness
theorem holds. This result is derived from the Paley-Wiener theorem \cite{paley}
and is the main tool in the study of Fourier transforms of weighted algebras.

Fix one more notation. For $x>0$ denote
$$
\log^+x=\max\{0,\log x\},\quad \log^-x=\min\{0,\log x\}.
$$
A weight $w$ on the real line is called non-quasianalytic, if
\beq\label{paley}
\int_{-\infty}^\infty {\log^+ w(t)\over 1+t^2} dt < \infty,
\eeq
and quasianalytic if this integral diverges.

Beurling \cite{beurl-1949} proved that \eqref{paley} is equivalent to regularity
of a weighted algebra $\L_1^w(\R)$. This result was further extended by Domar
\cite{domar-1956} to the case of any abelian locally compact group,
the definition of quasianalyticity slightly changed. According to Domar,
a weight $w$ on a group $G$ is called non-quasianalytic if for any $x\in G$
\beq\label{domar}
\sum_{n=1}^\infty {\log^+ w(nx)\over n^2} <\infty.
\eeq

\btm[{\cite[th. 2.11]{domar-1956}}]
Regularity of an algebra $\L_1^w(G)$ with a weight $w\ge1$.
is equivalent to the inequality \eqref{domar}.
\etm

The Domar's proof is easily extended to the case $p>1$, if one supposes
that the algebra $\lpw$ is translation invariant (theorem \ref{regul}).
For non-invariant algebras \eqref{domar} is not a criterion, see example
\ref{countex-tauber}. This distinction does not occur in the case $p=1$
as algebras $\L_1^w(G)$ are always translation invariant \cite[th. 3.1]{kuz-mz}.

Note that if we multiply the weight of an algebra $\lpw$ by a real-valued
character, $v=\chi w$, then we get an algebra $\L_p^v(G)$ which is isometrically
isomorphic to the original one. In particular, both algebras are regular
or not simultaneously. But the condition \eqref{domar} may cease to hold
after such a transition (consider, e.g., weights $w(t)=1+t^2$ and $v(t)=e^t(1+t^2)$
on the real line). It is known \cite[th.\ 3]{kuz-faa} that it is always
possible to find a character $\chi$ such that the new algebra $\L_p^v(G)$
with the weight $v=\chi w$ is contained in $\L_1(G)$. The latter inclusion
simplifies various details, and particularly, the condition \eqref{domar}
is equivalent to regularity only for algebras satisfying $\lpw\subset \L_1(G)$.
Thus, it should be kept in mind that the inequality \eqref{domar} is to be
tested after previous renorming of the weight.

Proof of theorem \ref{regul} requires the following lemma.

\blm\label{inf_w}
If $\lpw\subset \L_1(G)$ is an invariant algebra, then
$${\rm ess}\inf_G w>0.$$
\elm

\bpr
By theorem \cite[2.7]{feicht} $w$ may be chosen continuous, and by
proposition \cite[1.16]{edw} the function $L$ defined in \eqref{ledw}
is bounded on every compact set. Let $D=D^{-1}$ be a compact set of positive
measure, and $N$ such a number that $L_r\le N$ for $r\in D$.
Then for any $s\in G$, $r\in D$
\beq\label{wD}
{w(s)\over N}\le w(s+r) \le Nw(s).
\eeq
Suppose now that $\inf_G w=0$, i.e. $w(t_n)=x_n\to0$ for some sequence $t_n\in G$.
Then $w(t)\le Nx_n$ for $t\in t_n+D$, so that
$$
\int_G w^{-q}(t)dt \ge \int_{t_n+D} w^{-q}(t)dt \ge N^{-q} x_n^{-q}\mu(D) \to +\infty.
$$
But by assumption $\lpw\subset \L_1(G)$, what is by \cite[prop. 2]{kuz-faa}
equivalent to the inclusion $w^{-q}\in\L_1(G)$. We come to a contradiction.
\epr

We need also

\btm[{\cite[th. 3.3]{kuz-mz}}]\label{lp_invar}
Let $G$ be an abelian locally compact group, and let $\L_p^w(G)$ be an invariant algebra.
Then $w$ is equivalent to a submultiplicative function.
\etm

\bcor\label{p_submult}
Let $G$ be an abelian locally compact group, and let $\L_p^w(G)\subset \L_1(G)$ be
an invariant algebra. Then $w$ is equivalent to a submultiplicative function $v\ge1$.
\ecor

\bpr
Let $v=Cw$ be equivalent to $w$ and submultiplicative (th. \ref{lp_invar}).
By lemma \ref{inf_w} $\delta=\inf_G v>0$. If $\delta\ge1$, then $v$ is the desired
function. If $\delta<1$, then it is $v/\delta$.
\epr

\btm\label{regul}
Let $\lpw\subset \L_1(G)$ be an invariant algebra.
Regularity of this algebra is equivalent to the condition \eqref{domar}.
\etm

\bpr
{\it Sufficiency.}
By lemma \ref{p_submult} we may assume that $w\ge1$ and  $w$ is submultiplicative.
If now \eqref{domar} holds, then $\L_1^w(G)$ is a regular algebra,
and since $\lpw$ is a module over $\L_1^w(G)$ \cite[th.~1.1]{kuz-mz},
this algebra is also regular.

{\it Necessity.} Suppose the contrary, i.e. that the series \eqref{domar}
diverges for some $x\in G$. Let us assume the weight is continuous
(corollary \ref{p_submult}). The closed subgroup $G_x$ generated by $x$
is either compact or discrete in $G$ and therefore isomorphic to $\Z$ \cite[th.~9.1]{HR}.
In the first case series \eqref{domar} converges irrespective of regularity.
Suppose therefore that $G_x$ is discrete in $G$. Its dual group
$\hat G_x=\hat G/G_x^\perp$ may be identified with the unit circle, parameter for
which we take in $[-\pi,\pi]$. Pick a number $\e\in(0,\pi)$. By assumption
for any neighborhood of zero $U\subset \hat G$ there exists $f_0\in\lpw$
such that its support $\supp \hat f_0$ is compact and contained in $U$;
we take $f_0$ such that $U=(-\e,\e)+G_x^\perp$.

Take now any $f_1\in\lpw$ with compact support provided that
$\hat f_1\hat f_0\ne0$. Convolution $f=f_0*f_1$ may be estimated in an arbitrary
point $s$ using submultiplicativity of the weight:
\begin{gather*}%HR 20.10
|f(s)| \le \int \left|f_0(t)f_1(s-t){w(t)w(s-t)\over w(s)}\right|dt
 \le\\\le
  {1\over w(s)} \|f_0w\|_p \|f_1w\|_q = {C\over w(s)},
\end{gather*}
where the constant $C$ depends of course on the choice of $f_0$ and $f_1$.
Note that $\supp \hat f\subset (-\e,\e)+G_x^\perp$.
Restriction of $f$ onto subgroup $G_x$, $\phi_n=f(nx)$, satisfies
\beq\label{phi}
|\phi_n|\le {C\over w(nx)}, \quad \supp\hat\phi\subset (-\e,\e).
\eeq
Indeed, since $\hat\phi=T_{G_x^\perp}\hat f$
\cite[2.7.3]{rudin} (here $T$ denotes the averaging operator
\beq\label{mean}
T_N f(x)=\int_N f_x(t)dt = \int_N f(xt)dt,
\eeq
see \cite[\S 9]{weyl}),
then $\supp\hat\phi =\supp T_{G_x^\perp}\hat f\subset (-\e,\e)$.

Obviously $\phi$ can be extended to an entire function of exponential type
so that \cite[th. X, XII]{paley}
$$
\int_{-\infty}^\infty {\log^-|\phi(t)|\over 1+t^2}\,dt>-\infty.
$$
Submultiplicativity of $w$ will allow to conclude that the series \eqref{domar} converges.

Choosing $\phi_1\in \ell_1(\Z)$ so that the support $\supp\hat \phi_1$ is sufficiently
small, we may achieve that for $\tilde \phi = \phi\cdot \phi_1$ the condition \eqref{phi}
still holds, moreover, $\tilde\phi\cdot w\in\ell_1(\Z)$.
We may assume therefore that $\phi\cdot w\in\ell_1(\Z)$.

We show first that $\phi$ as a function on the real line changes little at small translations.
Let $|s|\le1$, $n\in\Z$. Denote $\xi_s(t)=e^{ist}$, then
$$
\phi(n+s)=(\xi_s\cdot\hat\phi)\,\check{\,}(n)=(\hat\xi_s*\phi)(n)\le
 \sum_k |\hat\xi_s(k)\phi(n-k)|.
$$
Due to \eqref{phi}
$$
|\phi(n-k)|\le {C\over w(nx-kx)} \le {Cw(kx)w(nx-kx)\over w(nx-kx)w(nx)}
 = {Cw(kx)\over w(nx)}.
$$
Convolution with $\xi_s$ is a multiple of $\xi_s$:
$$
(\xi_s*\hat\phi)(t) = {1\over2\pi}\int_{-\pi}^\pi \hat\phi(r)e^{i(t-r)s}dr
 = \xi_s(t)\phi(s),
$$
i.e. $\xi_s = \xi_s*\hat\phi/\phi(s)$ и $\hat\xi_s = \hat\xi_s\cdot\phi/\phi(s)$.
Thus,
$$
|\phi(n+s)| \le  \sum_k \Big|\hat\xi_s(k){\phi(k)\over\phi(s)}\Big|\cdot{Cw(kx)\over w(nx)}
\le {C_1\over w(nx)|\phi(s)|} \|\phi\cdot w\|_1 \le {C_2\over w(nx)}
$$
for $s$ in sufficiently small interval $(-\d,\d)$, as $\phi(0)=w(0)\ne0$
and $\phi$ is continuous as a function on the real line.

But now
\begin{gather*}
\sum_{n=1}^\infty {\log^+ w(nx)\over n^2} \le
{1\over2\d}\, \sum_{n=1}^\infty \,\int_{n-\d}^{n+\d} {\log^+ C_2-\log^- |\phi(s)|\over n^2}ds
\le\\
\le C_3-C_4\int_{-\infty}^\infty {\log^- |\phi(t)|\over 1+t^2}dt
<+\infty.
\end{gather*}
\epr

On the real line regularity of an algebra $\L_p^w(\R)$ is equivalent to convergence
of the integral \eqref{paley} even if the weight is not submultiplicative.
A simple proof of sufficiency in the case $p=1$ is presented in a paper of
E.~A.~Gorin \cite{gorin} and may be literally repeated for $p>1$.
Necessity may be proved, for example, in the following way. Let $f\in\L_p^w(\R)$
be such that the support $\supp\hat f$ is compact. Then by Paley-Wiener theorem
$$
\int_{-\infty}^\infty {\log^-|f(t)|\over 1+t^2}\,dt>-\infty.
$$
By definition of the weighted space $fw=\phi\in\L_p(\R)$. Since
$\log^+|\phi|\le|\phi|$, then $\log^+|\phi|\in\L_p(\R)$. From the other side,
$(1+t^2)^{-1}\in\L_q(\R)$ for all $q\ge1$. Thus, from the equality $w=\phi/f$ we get:
$$
\int_{-\infty}^\infty {\log^+w(t)\over 1+t^2}\,dt      \le
\int_{-\infty}^\infty {\log^+|\phi(t)|\over 1+t^2}\,dt -
\int_{-\infty}^\infty {\log^-|f(t)|\over 1+t^2}\,dt    <+\infty.
$$

Submultiplicativity of the weight is essential when we pass from the integral
\eqref{paley} to the series \eqref{domar}, and it would be natural to expect
that in its absence condition \eqref{domar} is no longer equivalent to
regularity. It is indeed so what shows the following example.

\bex\label{countex-tauber}
Consider the unit circle $\T$ with parameter $t\in[0,1)$ and the weight $w(t)=t^{1/4}$
on it. The space $\L_2^w(\T)$ is an algebra because condition \eqref{eqw} holds.
Algebra $\L_2^w(\T)$ is regular because it contains all exponents of type
$f_n(t)=e^{int}$ with Fourier transforms $\hat f_n=\delta_n$.
But the condition \eqref{domar} does not hold. In order to show it, choose
a number $\a\in(0,1)$ with good rational approximations, e.g.
$\a=\sum_{n=1}^\infty q_n^{-1}$, where $q_1=2$, $q_n> 2q_{n-1}\exp(q_{n-1}^2)$,
and all $q_n$ are integer. Then $\{q_n\a\} < 2q_n/q_{n+1} < e^{-q_n^2}$. As group
operation on $[0,1)$ is the fractional part of ordinary sum,
$$
\sum_{n=1}^\infty {|\log w(n\a)|\over n^2}
\ge \sum_{n=1}^\infty {|\log w(q_n\a)|\over q_n^2}
= \sum_{n=1}^\infty {|\log \{q_n\a\}|\over 4q_n^2}
 \ge \sum_{n=1}^\infty {|-q_n^2|\over 4q_n^2} = +\infty.
$$
\eex

From results of Domar \cite[1.5]{domar-1956} follows
\btm\label{spec_reg}
The spectrum of a regular algebra $\L_p^w(G)\subset \L_1(G)$ coincides with the
dual group $\hat G$.
\etm

\bcor
Spectrum of a regular algebra $\L_p^w(G)$ is ho\-meo\-morph\-ic to the dual group $\hat G$.
\ecor

\bpr
It follows from the previous theorem and theorem \cite[th.\ 3]{kuz-faa}.
\epr

Thus, a typical example of a weight of an irregular algebra on the real line
is given by $w(t)=(1+t^2)e^{|\,t|}$. Spectrum in this example is the strip
$-1\le {\mathfrak Re} z\le 1$ (a proof for $p=1$, valid also for submultiplicative
weight and $p>1$, see in \cite[\S 18]{grsh}). Next example shows that the
spectrum of an irregular algebra can be also equal to the dual group.

\bex
Let $w(t)=(1+t^2)\exp\big(|\,t|/\log(e+|\,t|)\big)$. Then series \eqref{domar}
diverges for any $t\ne0$, as $\log w(t) = O(\log t) + |\,t|/\log(e+|\,t|)$.
At the same time $\L_2^w(\R)$ is an algebra (inequality \eqref{eqw} is checked straitforwardly).
Since $\lim_{t\to\pm\infty} t^{-1} \log w(t)=0$,
spectrum of this algebra is the imaginary axis \cite[\S 18]{grsh},
which may be identified with the dual group $\R$.
\eex

\section{Construction of regular algebras}\label{constr-reg}

In this section we show that regular weighted algebras with $p>1$ may be
constructed on every $\s$-compact abelian group. Since
$\s$-compactness is (in abelian case) necessary for the existence of weighted
algebras with $p>1$ (th. \cite[1.1]{kuz-mz}), we see that if weighted algebras
do exist on a given group, then there are also regular ones among them.
In the case $p=1$ there is no problem because the classical algebra $\L_1(G)$ is regular.

We say that the weight $w$ grows polynomially if there exists $d\in\N$ such
that for all $x\in G$
$$
w(nx)=O(n^d),\quad n\to\infty.
$$
It is clear that for such a weight Domar's condition \eqref{domar} holds.
On a compactly generated abelian group it makes no trouble to construct such
a weight. The task becomes more complicated when the group is not compactly generated,
e.g. in the case of rationals $\Q$ or $p$-rationals $\Z(p^\infty)$.
Separate lemmas \ref{Q} and \ref{fin} are devoted to these two groups.
The group $\Z(p^\infty)$ for prime $p$ is defined as the set
$\{k/p^n: k\in\Z,n\in\Z_+\}$ with addition modulo 1.
In theorem \ref{exist-abel} we show with the help of structure theory
that a polynomially growing weight may be constructed on any $\s$-compact
abelian group.

We construct the weight via auxiliary function $u=w^{-q}$ so that \eqref{eqw} holds.
This guarantees that all spaces $\lpw$ are convolution algebras.
We list properties that will be required from $u$ in the following lemmas.
This function must be:
\bit
\item[(a)] Positive: $u>0$
\item[(b)] $u * u \le u$
\item[(c)] Even: $u(x)=u(-x)$
\item[(d)] Decrease polynomially: there exists $d$ such that
 $1/u(nx)=O(n^d),\quad n\to\infty$ for all $x$.
\eit

\blm\label{fin}
Let $G$ be equal to the union of its nested
finite subgroups $G_n$, $n\in\N$: $G=\cup G_n$,
$G_n\subset G_{n+1}$ for all $n$.
Then on $G$ exists a weight with properties (a)--(d).
\elm

\bpr Let $|M|$ denote cardinality (which is also Haar mea\-su\-re in this case)
of a set $M\subset G$.
We may construct the weight with the help of any decreasing sequence
$\phi_n>0$ provided that $\phi=\sum \phi_n|G_n|<\infty$. Denote
$U_j=G_j\setminus G_{j-1}$, assuming $G_0=\emptyset$. Then
$G=\cup U_j$, and we define $u$ by $u|_{U_n} = \phi_n$.

Properties (a), (c) are obvious. Since $1/u(nx)$ is bounded for any $x$,
the decrease condition (d) is satisfied. Finally, we show that (b) holds.
It is clear that
$$
(u*u)(x) = \sum_{j=1}^\infty \sum_{y\in U_j} \phi_j u(x-y).
$$
Let $x\in U_n$, $y\in U_j$. If $j<n$ then $x-y\in U_n$; if
$j>n$, then $x-y\in U_j$. The set $U_n$ splits into disjoint
union of the sets $x+U_j$, $j<n$, and their complement
$U_x$, for which $x-U_x\subset U_n$.
Since $u$ is constant on all these sets, we get that
\begin{gather*}
(u*u)(x)\! = \!\sum_{j<n}\! |U_j|\phi_j \phi_n +\! \sum_{j<n}\!|x+U_j|\phi_n\phi_j
 + |U_x|\phi_n^2 +\! \sum_{j>n}\! |U_j|\phi_j^2
\le\\\le
 2\!\sum_{j=1}^\infty\! |U_j|\phi_j\phi_n \le \!2\phi_n\phi\! =\! 2u(x)\phi,
\end{gather*}
i.\ e.\ (b) holds after transition from $u$ to $2\phi u$.
\epr

\bcor\label{Zp}
On the group $\Zp$ for any prime $p$ there exists a weight with properties (a)--(d).
\ecor

\blm\label{Q}
On the group $\Q$ there exists a weight with properties (a)--(d).
\elm

\bpr Represent $\Q$ as a union of nested subgroups $Q_n =
\Z/t_n$, $n\in \N$ --- one may take, e.g., $t_n=n!$. Denote
$U_n=Q_n\setminus Q_{n-1}$, assuming $Q_0=\emptyset$. Let us use functions
on $\Z$
\beq
\bar n=\max\{1,|n|\} \text{ and }\s(n)=1/\bar n^2.
\eeq
Let $C_2$ be such that
\beq\label{n^a}
\sum_{n=-\infty}^\infty {1\over \s(n)\s(m-n)} \le {C_2\over\s(m)}
\eeq
(its existence is easy to verify).

Pick now a decreasing sequence $\phi_n>0$ such that
$\sum \phi_n t_n=\phi<\infty$. For $q\in U_n$ we put
$u(q) = \phi_n \sigma(\lfloor q\rfloor)$, where $\lfloor x\rfloor$
is the even extension of $[x]$ from the real half-line.

It is obvious that this weight is even and positive.
Any $q\in\Q$ belongs to some subgroup $Q_m$, so that
$u(nq)\ge \phi_m\s(\lf nq\rf)=\phi_m/\lf nq\rf^2$. Obviously this sequence
decreases at most quadratically on $n$.

Similarly to the proof of the previous lemma $(u*u)(q)$ may be split into the sum
$$(u*u)(q) = \sum_{j<n} (S_j + S'_j) + S_n + \sum_{j>n} S_j,$$
where for $j<n$
$$
S_j = \sum_{r\in U_j} \phi_j \phi_n \sigma(\lf r\rf)\sigma(\lf q-r\rf)
 = S'_j = \sum_{r\in q+U_j}\phi_n\phi_j \sigma(\lf r\rf)\sigma(\lf q-r\rf),
$$
for $j>n$
$$
S_j = \sum_{r\in U_j}\phi_j^2 \sigma(\lf r\rf)\sigma(\lf q-r\rf)
$$
and, finally,
$$
S_n = \sum_{r\in U_n\setminus \mathop{\cup}\limits_{j<n}(q+U_j)}
 \phi_n^2\sigma(\lf r\rf)\sigma(\lf q-r\rf).
$$
We estimate now these sums separately. For $j<n$
$$
{S_j\over u(q)} = \phi_j \phi_n \sum_{k=0}^\infty
 \sum_{r\in U_j\cap\Big((-k-1,-k]\cup[k,k+1)\Big)}
  {\sigma(\lf r\rf)\sigma(\lf q-r\rf)\over\phi_n\sigma(\lf q\rf)}.
$$
If $\lf r\rf=k$, then $\lf q-r\rf$ lies between $\lf q\rf-k-1$ and
$\lf q\rf-k+1$. A rough estimate $\sigma(l\pm1)\le 4\sigma(l)$ holds.
Note also that $|U_j\cap(-k-1,-k]|=|U_j\cap[k,k+1)|\le t_j$.
Moreover, $\lf x\rf$ equals 0 on two integer segments. Thus (recall that
the notation $C_2$ is introduced in \eqref{n^a}),
$$
{S_j\over u(q)} \le 8\phi_j t_j
 \sum_{k\in\Z}{\sigma(k)\sigma(\lf q\rf-k)\over \sigma(\lf q\rf)}
 \le 8C_2\phi_j t_j\equiv C\phi_j t_j.
$$
Similarly, for $j\ge n$ we get $S_j\le u(q)\cdot C t_j\phi_j^2/\phi_n
\le u(q)\cdot C t_j\phi_j$ (the latter is true due to decrease of $\phi_j$).
So,
$$
(u*u)(q) \le u(q)\cdot2\sum_{j=1}^\infty C\phi_j t_j = u(q)\cdot 2\,C\phi,
$$
i.\ e.\ (b) holds after transition to $2\,C\phi \,u$.
\epr

A direct sum $\oplus_{\a\in A}G_\a$ of groups $G_\a$ is usually defined
\cite[B7]{rudin} as the subgroup in their direct product
$\prod_{\a\in A}G_\a$ consisting of all elements with finitely many
nonzero coordinates. Next lemma constructs a weight on the countable direct
sum by weights on the summands. In fact, the proof remains true even
in non-commutative case.

\blm\label{weak_prod}
Let $G_j$, $j\in\N$ be discrete groups with weights $u_j$ satisfying (a)--(d).
On the discrete direct sum $G=\oplus_{j\in\N} G_j$ there exists a weight
$u$ satisfying the same conditions.
\elm

\bpr
An element $x\in G$ is defined by its coordinates $x_j\in G_j$.
We denote $s(x)=\{j\in\N: x_j\ne0\}$ and write also $s_x$
instead of $s(x)$. This set is finite for any $x$.
Denote also $G_j^\times = G_j\setminus\{0\}$.
In order to define the weight we need a sequence $\a_j>0$
and a function $a:s\to a_s$ on the set $\FN$ of all finite subsets of the
set of natural numbers; $\a_j$ and $a_s$ will be specified later.
Now, define $u$ as
\beq\label{udef_abel}
u(x) = a_{s(x)} \prod_{j\in s(x)}\a_ju_j(x_j).
\eeq
It is obviously even. We will choose $a_s$ positive, so that $u>0$.
Note also that $s(nx)=s(x)$ for all $x\in G$, and $j$-th coordinate of $nx$ is
$nx_j$. Thus, decrease condition (c) for $u$ follows trivially from
the same condition holding for every $u_j$.

For (b), the following inequalities are sufficient.
Numbers $\a_j$ should be small:
\begin{align}
\forall j\in\N\quad 0<\a_j&<1,            \label{alphaone_abel}\\\sk{2pt}
\prod_{j=1}^\infty (1+\a_j^2\cdot (u_j*u_j)(0))  &<2, \qquad     \label{alphaujnul}\\
\qquad
\prod_{j=1}^\infty (1+\a_j)         &<2.   \label{alphaKj}
\end{align}
It is obvious that such $\a_j$ exist. From $a$ we require that for all $s\in\FN$
\begin{align}
0< a_s &\le1,                           \label{aone}\\
\sk{2pt}
\kern40pt a_{s\cup v} &\le a_s  \quad\forall v\in\FN,    \label{aunion}\\
\sum_{v\subset s}\displaystyle {a_va_{s\setminus v}\over a_s} &\le{1\over4}. \label{asubset}
\end{align}

One may take $a_s=\e_1/(\sum_{j\in s}j!)$ with $\e_1=e^{-2}/8$.
For the empty set we put $a_\emptyset=\e_1$.
Properties \eqref{aone}, \eqref{aunion} hold obviously. Check now \eqref{asubset}.
Note that for nonempty $v\subsetneq s$ we have
$a_s^{-1} = a_v^{-1}+a_{s\setminus v}^{-1}$. Therefore
$$
\sum_{v\subset s} {a_va_{s\setminus v}\over a_s}
 = 2a_\emptyset + \sum_{v\subsetneqq s} a_va_{s\setminus v}
  \Big({1\over a_v}+{1\over a_{s\setminus v}}\Big)
 = 2\sum_{v\subsetneqq s} a_v.
$$
Further, $a_v=\e_1/(\sum_{j\in v}j!) \le \e_1/(\sup v)!$.
assuming $\sup \emptyset=0$. Thus,
$$
\sum_{v\subsetneqq s} a_v \le \sum_{v\subsetneqq s} {\e_1\over(\sup v)!}
 \le \sum_{k=0}^{\sup s} {\e_1 2^k\over k!}
 \le \e_1 e^2 = {1\over8},
$$
whence \eqref{asubset} follows.

Estimate now the convolution $u*u$.
Take $x\in G$. Every $x'\in G$ may be represented then in the form
$x'=y+z$, where $s_y\subset s_x$, $s_z\subset \N\setminus s_x$.
According to the definition \eqref{udef_abel},
$$
u\/(x')=u\/(y+z) = u\/(y) {a_{s_y\cup s_z}\over a_{s_y} } \prod_{j\in s_z}\a_j u_j(z_j)
\overset{\textstyle\eqref{aunion}}\le\,
  u(y) \prod_{j\in s_z}\a_j u_j(z_j).
$$
As $s_{x-y-z}=s_{x-y}\cup s_z$ and
$(x-y-z)_j=-z_j$ for $j\in s_z$, using the fact that all weights are even,
\begin{align*}
u\/(x-x')=u\/(x-y-z)
= u\/&(x-y) {a_{s_{x-y}\cup s_z}\over a_{s_{x-y}} }  \prod_{j\in s_z}\a_j u_j(z_j)
  \le\\
\overset{\textstyle\eqref{aunion}}\le\,
  u&(x-y) \prod_{j\in s_z}\a_j u_j(z_j).
\end{align*}
Consider subgroup $Q_x=\prod_{j\in s_x} G_j$ extracted by the support of $x$.
An element $y$ belongs to this subgroup iff $s_y\subset s_x$.
Using the above inequalities, we get:
\begin{align*}
u*u\,(x) = \sum_{x'\in G} u\/(x')u\/(x-x') \le
 \sum_{y\in Q_x} u\/(y)u\/(x-y) \kern-5pt
  \sum_{z\in G \atop \scriptstyle s_z\subset \N\setminus s_x}
   \prod_{j\in s_z}\a_j^2 u_j^2(z_j)
\le\\\le
 \sum_{y\in Q_x} u\/(y)u\/(x-y)
   \prod_{j=1}^\infty\big(1+\a_j^2\cdot (u_j*u_j)(0)\big)
  \overset{\textstyle\eqref{alphaujnul}}<\,
   2\sum_{y\in Q_x} u\/(y)u\/(x-y).
\end{align*}

It remains thus to estimate sum over $Q_x$. We will write further $s$ instead of $s_x$.
Let $y\in Q_x$ and $s_y=v$. The difference $x-y$ has some support $s_{x-y}$. Denote
$r=v\cap s_{x-y}$. Then %для $j\in s\setminus r$  $y$ и $x\!-\!y$
coordinates off $r$ are determined uniquely:
$$\begin{cases}
y_j=x_j,\ \;(x\!-\!y)_j=0 \qquad&\text{ for } j\in v\setminus r,\\
y_j=0,\ \;(x\!-\!y)_j=-x_j \qquad&\text{ for } j\in s\setminus v.
\end{cases}$$
Taking into account all weights are even, we note that in the following fraction
after cancellation it remains only
\begin{align*}
{u\/(y)u\/(x-y) \over u\/(x)} &=
 {a_v a_{r\cup (s\!\setminus\! v)} \over a_s}
 \prod_{j\in r} \a_j {u_j(y_j)u_j(x_j-y_j) \over u_j(x_j)}
\le\\  &\overset{\textstyle\eqref{aunion}}\le\;\;
 {a_v a_{s\!\setminus\! v} \over a_s}
 \prod_{j\in r} \a_j {u_j(y_j)u_j(x_j-y_j) \over u_j(x_j)}.
\end{align*}
Estimate sum over $Q_x$:
$$\displaylines{
\sum_{y\in Q_x} {u\/(y)u\/(x-y)\over u\/(x)} =
 \sum_{v\subset s} \sum_{s(y)=v} {u\/(y)u\/(x-y)\over u\/(x)} =
\cr=
\sum_{v\subset s} {a_v a_{s\!\setminus\! v} \over a_s}
 \sum_{r\subset v}
  \sum_{s(y)=r \atop y_j\ne x_j}
   \prod_{j\in r} \a_j{u_j(y_j)u_j(x_j-y_j)\over u_j(x_j)}
=\cr=
\sum_{v\subset s} {a_v a_{s\!\setminus\! v}\over a_s}
 \sum_{r\subset v}
 \prod_{j\in r}\a_j \sum_{0\ne y_j\ne x_j} {u_j(y_j)u_j(x_j-y_j)\over u_j(x_j)}
\overset{\textstyle(b)}\le\,
\sum_{v\subset s}
 {a_v a_{s\!\setminus\! v}\over a_s}
 \sum_{r\subset v} \prod_{j\in r}\a_j
=\cr=
\sum_{v\subset s}
 {a_v a_{s\!\setminus\! v}\over a_s}
  \prod_{j\in v} (1+\a_j)
\overset{\textstyle\eqref{alphaKj}}\le\,
2\sum_{v\subset s}
 {a_v a_{s\!\setminus\! v}\over a_s}
\overset{\textstyle\eqref{asubset}}\le\,
 {1\over2}.
}$$

Thus, for all $x\in G$
$$
u*u\,(x) \le u(x),
$$
i.e. (b) holds.% End of the proof.
\epr

We can now construct a weight on any countable group, using the structure theorem.

\blm\label{exist-discr}
On every countable discrete abelian group $G$ exists a weight with properties (a)--(d).
\elm

\bpr
It is known that $G$ may be embedded as a subgroup into a divisible group
$H$, which is also countable. By structure theorem for abelian groups
\cite[suppl.\ A]{HR} $H$ is isomorphic to the direct sum of copies of rationals
$\Q$ and $p$-rational numbers $\Z(p^\infty)$ with prime $p$.
Since $H$ is countable, the number of summands is countable.
By lemmas \ref{Q}, \ref{Zp} on all summands exist weights with properties (a)--(d).
By lemma \ref{weak_prod} the same is true for $H$. If we restrict now the weight
from $H$ onto $G\subset H$, properties (a)--(d) with remain true.
\epr

\btm\label{exist-abel}
On every $\sigma$-compact locally compact abelian group $G$ for every $p>1$
there exists an invariant algebra $\lpw$ with polynomially growing weight.
\etm

\bpr By structure theorem (\cite[24.30]{HR}) $G$ is topologically isomorphic
to $\R^d\times H$, where $d$ is a nonnegative integer and $H$ contains
a compact open subgroup $E$. The quotient group $H/E$ is discrete, and
by $\sigma$-compactness of $G$ it is countable.

Let functions $u_R$ on $\R^d$ and $u_H$ on $H/E$ satisfy (a)--(d).
Existence of $u_H$ is proved in lemma \ref{exist-discr}, and $u_R$ may be taken equal to
$$
u_R(x)={1\over (1+x_1^2)\cdots(1+x_d^2)}.
$$
The function $u$ on $G$,
$$
u(r,h) = u_R(r) u_H(h+E),
$$
also satisfies (a)--(d).
But now for any $p>1$ the space $\lpw$ with the weight $w=u^{-q}$ is an algebra,
its weight growing polynomially. Moreover, $u_H^{-q}$ is submultiplicative
because it defines an algebra \cite[remark~3.1]{kuz-mz}. We see that $u^{-q}$ and $w$
are also submultiplicative, thus the algebra $\lpw$ is invariant.
\epr

Now it follows

\btm\label{reg-abel}
On every $\sigma$-compact locally compact abelian group for any $p>1$
there exists a regular algebra $\lpw$.
\etm

Author thanks deeply E. A. Gorin for the statement of the problem and
useful discussions.


\begin{thebibliography}{99}

\bibitem{beurl-1949}
Beurling A.
On the spectral synthesis of bounded functions,
{\it Acta Math.}
{\bf 81} (1949), 225--238.

\bibitem{domar-1956}
Domar Y.
Harmonic analysis based on certain commutative algebras,
{\it Acta Math.}
{\bf 96}, no.~2 (1956), 1--66.

\bibitem{edw}
Edwards R.~E. The stability of weighted Lebesgue spaces.
{\it Trans. Amer. Math. Soc.}
{\bf 93} (1959), 369--394.

\bibitem{feicht}
Feichtinger H. G.
Gewichtsfunktionen auf lokalkompakten Gruppen,
{\it Sitzber. \"Osterr. Akad. Wiss. Abt. II,}
{\bf 188}, no.~8--10 (1979), 451--471.

\bibitem{grsh}
%Гельфанд И.~М., Райков Д.~А., Шилов Г.~Е.
%{\it Коммутативные нормированные кольца.} М.: Физматгиз, 1960.
Gel'fand I. M., Ra\u\i kov D. A., \v Silov G. E.
{\it Kommutativnye normirovannye kol'tsa [Commutative normed rings]}.
Fiz.-Mat. Lit., Moscow, 1960.

\bibitem{gorin}
%Горин Е. А. Исследования Г. Е. Шилова по теории коммутативных банаховых алгебр
%и их дальнейшее развитие. Успехи мат. наук {\bf 33}, №4 (202), 169--188 (1978).
Gorin E. A. G. E. \v Silov's investigations in the theory of commutative Banach algebras
and their subsequent development. {\it Uspekhi Mat. Nauk} {\bf 33} (1978), no. 4 (202),
169--188, 256.

\bibitem{HR}
Hewitt, E.; Ross, K.A.
{\it Abstract harmonic analysis I, II.}
Springer--Verlag, 3rd printing, 1997.

\bibitem{kuz-faa}
%Кузнецова Ю.~Н. Весовые $L_p$-алгебры на группах.
%{\it Функц. анализ и его прил.,} {\bf 40}, № 3 (2006), 82-85.
Kuznetsova Yu. N. Weighted $L_p$-algebras on groups. {\it Funkts. Anal. i Pril.} {\bf 40} (2006),
no. 3, 82--85; transl. in {\it Funct. Anal. Appl.} {\bf 40} (2006), no. 3, 234--236.

\bibitem{kuz-mz}
%Кузнецова Ю.~Н. Инвариантность весовых алгебр $L_p^w(G)$. {\it Матем. заметки} (в печати).
Kuznetsova Yu. N. Invariant weighted algebras $L_p^w(G)$. {\it Mat. zametki} (in print).

\bibitem{loom}
Loomis L. {\it An introduction to abstract harmonic analysis}.
D.~Van Nostrand Company, 1953.

\bibitem{paley}
Paley R., Wiener N. {\it Fourier transforms in the complex domain}.
Amer. Math. Soc., 1934.

\bibitem{rudin}
Rudin W.
{\it Fourier analysis on groups.}
NY: Interscience, 1962.

\bibitem{shilov}
%Шилов Г. Е. О регулярных нормированных кольцах.
%{\it Труды Матем. ин-та им. Стеклова} {\bf 21} (1947), 1--118.
Shilov G. E. On regular normed rings. {\it Trav. Inst. Math. Stekloff} {\bf 21} (1947).

\bibitem{weyl}
Weil A. {\it L'Int\'egration dans les groupes topologiques
et ses applications}. Paris: Hermann, 1940.

\end{thebibliography}
\end{document}